\def\draft{n}
\documentclass{amsart}
\usepackage[headings]{fullpage}
\usepackage{amssymb,epic,eepic,epsfig,amsbsy,amsmath}

\theoremstyle{plain}

\newtheorem{theorem}{Theorem}
\newtheorem{proposition}{Proposition}[section]
\newtheorem{lemma}[proposition]{Lemma}
\newtheorem{corollary}[proposition]{Corollary}

\newtheorem{conjecture}{Conjecture}

\theoremstyle{definition}
\newtheorem{definition}[proposition]{Definition}

\newtheorem{question}{Question}

\theoremstyle{remark}

\newtheorem{remark}[proposition]{Remark}

\def\printname#1{
    \if\draft y
        \smash{\makebox[0pt]{\hspace{-0.5in}
            \raisebox{8pt}{\tt\tiny #1}}}
    \fi
}

\newcommand{\psdraw}[2]
         {\begin{array}{c} \hspace{-1.3mm}
    \raisebox{-4pt}{\epsfig{figure=draws/#1.eps,width=#2}}
    \hspace{-1.9mm}\end{array}}

\newlength{\standardunitlength}
\catcode`\@=11
\long\def\@makecaption#1#2{%
     \vskip 10pt

\setbox\@tempboxa\hbox{
       \small\sf{\bfcaptionfont #1. }\ignorespaces #2}%
     \ifdim \wd\@tempboxa >\captionwidth {%
         \rightskip=\@captionmargin\leftskip=\@captionmargin
         \unhbox\@tempboxa\par}%
       \else
         \hbox to\hsize{\hfil\box\@tempboxa\hfil}%
     \fi}
\font\bfcaptionfont=cmssbx10 scaled \magstephalf
\newdimen\@captionmargin\@captionmargin=2\parindent
\newdimen\captionwidth\captionwidth=\hsize
\catcode`\@=12

\newcommand{\qbinom}[2]{\text{$\left[\begin{array}{c}#1\\ #2\end{array}
\right]$}}

\def\lbl#1{\label{#1}\printname{#1}}

\def\BN{\mathbb N}
\def\BZ{\mathbb Z}

\def\BQ{\mathbb Q}

\def\BC{\mathbb C}

\def\A{\mathcal A}
\def\cA{\mathcal A}
\def\cR{\mathcal R}

\def\cL{\mathcal L}

\def\W{\mathcal W}

\def\calI{\mathcal I}

\def\La{\Lambda}

\def\la{\langle}
\def\ra{\rangle}

\def\w{\omega}

\def\d{\delta}

\def\ti{\widetilde}

\def\longto{\longrightarrow}

\def\Z{\BZ}
\def\w{\omega}

\def\msl{\mathfrak{sl}}

\def\ihs{integer homology 3-sphere}
\def\Lhat{\widehat{\La}}

\def\Om{\Omega}
\def\ev{\mathrm{ev}}

\def\Zp{Z^{\mathrm{pert}}}
\def\SL{\mathrm{SL}}
\def\fsl{\mathfrak{sl}}

\newcommand{\LL}{\mathcal L}

\newcommand{\Ha}{\widehat{\Z[q]}}
\newcommand{\RR} {\mathbf R}

\newcommand{\Lp}{\mathbf L}

\begin{document}


\title[Is the Jones polynomial of a knot really a polynomial?]{
Is the Jones polynomial of a knot really a polynomial?}

\author{Stavros Garoufalidis}
\address{School of Mathematics \\
         Georgia Institute of Technology \\
         Atlanta, GA 30332-0160, USA \\
         {\tt http://www.math.gatech} \newline {\tt .edu/$\sim$stavros } }
\email{stavros@math.gatech.edu}
\author{Thang TQ L\^e}
\address{School of Mathematics \\
         Georgia Institute of Technology \\
         Atlanta, GA 30332-0160, USA \\
{\tt http://www.math.gatech} \newline {\tt .edu/$\sim$letu } }
\email{letu@math.gatech.edu}

\thanks{The authors were supported in part by NSF. \\
\newline
1991 {\em Mathematics Classification.} Primary 57N10. Secondary 57M25.
}

\date{This edition: January 5, 2006
\hspace{0.5cm} First edition: January 5, 2006.}

\dedicatory{Dedicated to Louis Kauffman on the occasion of his 60th birthday}

\begin{abstract}
The Jones polynomial of a knot in 3-space is a Laurent polynomial in $q$,
with integer coefficients. Many people have pondered why is this so,
and what is a proper generalization of the Jones polynomial for knots
in other closed 3-manifolds. Our paper centers around this question.
After reviewing several existing definitions of the Jones polynomial,
we show that the Jones polynomial is really an analytic function, in
the sense of Habiro. Using this, we extend the holonomicity properties
of the colored Jones function of a knot in 3-space to the case of a
knot in an integer homology sphere, and we formulate an analogue of
the AJ Conjecture.
Our main tools are various integrality properties of topological
quantum field theory invariants of links in 3-manifolds, manifested
in Habiro's work on the colored Jones function.
\end{abstract}

\maketitle

\tableofcontents


\section{Introduction}
\lbl{sec.intro}

\subsection{The Jones polynomial of a knot}
\lbl{sub.jones}

In 1985 Jones discovered the celebrated {\em Jones polynomial} of a
knot/link in 3-space, see \cite{J}. The framed version of the Jones
polynomial of a {\em framed} oriented knot/link may be uniquely
defined by the following {\em skein theory}:
$$
q^{1/4} J_{\psdraw{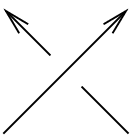}{.15in}}- q^{-1/4} J_{\psdraw{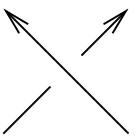}{.15in}}=
(q^{1/2}-q^{-1/2}) \, J_{\psdraw{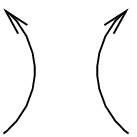}{.15in}}, \qquad
J_{\text{unknot}}=q^{1/2}+q^{-1/2}.
$$
From the above definition, the Jones polynomial $J_L$ of an oriented
link $L$ lies in $\BZ[q^{\pm 1/4}]$.

An important and trivial observation is that knots may be considered
to lie in other closed 3-manifolds other than $S^3$. The paper
is centered around the following question:
\begin{question}
\lbl{que.1}
What is the Jones polynomial of a knot $K$ in a closed 3-manifold $N$?
In particular, is the Jones polynomial of a knot in a closed 3-manifold
really a polynomial?
\end{question}
Before we reveal the answer, let us review some alternative
definitions of the Jones polynomial, namely the Kauffman bracket
approach, the quantum group approach, the TQFT
(Witten-Reshetikhin-Turaev) approach, and the perturbative TQFT
approach. Note that being polynomial is closely related to the
integrality discussed in \cite{Le3}.

\subsection{The Kauffman bracket of a knot}
\lbl{sub.kauffmanb}

Soon after Jones's discovery, Kauffman gave a reformulation of the Jones
polynomial in terms of the {\em Kauffman bracket}
$$
\la \cdot\ra : \text{Framed unoriented links in} \, S^3 \longto
\BZ[A^{\pm 1}]
$$
The Kauffman bracket is defined by the following skein theory (see \cite{Ka}):
$$
\left\la \psdraw{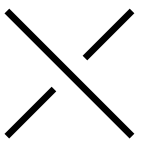}{0.5in} \right\ra =A \left\la \psdraw{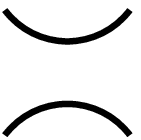}{0.5in}
\right\ra +A^{-1} \left\la \psdraw{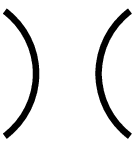}{0.5in} \right\ra, \qquad
\left\la \mathbf{\bigcirc} \right\ra = -(A^2+A^{-2})
$$
and the relation between the Kauffman bracket and the Jones
polynomial is the following: if $L$ is an oriented $m$-component
link projection with writhe $w(L)$, then
$$
J_L|_{q^{1/4}=A}=(-1)^{m}(-A^3)^{-w(L)} \la L \ra
$$
Since the Kauffman bracket takes values in $\BZ[A^{\pm 1}]$, it
implies once again that the Jones polynomial of an oriented
knot/link in 3-space takes values in $\BZ[q^{\pm 1/4}]$. Although
the Kauffman bracket can be generalized to framed unoriented links
in a closed 3-manifold, the resulting {\em Kauffman bracket skein
module} is an interesting and not well-understood object. Most
important, the Kauffman bracket skein module of a closed 3-manifold
is not a ring, but rather a module over the ring $\BZ[A^{\pm 1}]$.

\subsection{The quantum group approach}
\lbl{sub.qgroup}

Using quantum group technology (see \cite{RT,Tu}), one may define
the {\em colored Jones polynomial} of a framed oriented knot in
3-space. Let us postpone the technical details to a later Section
\ref{sub.qgroups} and concentrate on a main idea: the notion of {\em
color}.

Suppose $L$ is a framed $r$-component link in $S^3$ with ordered
components. The Jones polynomial of $L$ is a powerful
invariant that takes values in the Laurent polynomial ring
\begin{equation}
\lbl{eq.R}
\cR:=\BZ[q^{\pm 1/4}].
\end{equation}
An even more powerful invariant is the {\em colored Jones function}
$$
J_L :\BN^r \longto \cR
$$
which encodes the Jones polynomial of $L$ together with its
parallels. Here $\BN$ is the set of positive integers, and
$J_{L}(n_1,\dots,n_r)$ is the quantum $\fsl_2$ invariant of the link
$L$ with colors the irreducible $\fsl_2$-modules of dimensions
$n_1,\dots,n_r$. When  $n_1=\dots=n_r=2$, the polynomial
$J_{L}(n_1,\dots,n_r)$ is the Jones polynomial. Here we use the
normalization such that when $U$ is the unknot of 0 framing,
$$
J_U(n) = [n] := \frac{q^{n/2} - q^{-n/2}}{q^{1/2}- q^{-1/2}}.
$$

\subsection{The TQFT and the perturbative approach}
\lbl{sub.TQFT}

Let us fix a framed oriented knot $K$ in a 3-manifold $N$ and a
nonnegative integer $n \in \BN$. The Witten-Reshetikhin-Turaev (WRT)
invariant $Z_{(N,K),\xi,V_n} \in \BQ(\xi)$ is a generalization of
the colored Jones polynomial and can be defined under the framework
of {\em topological quantum field theory} (TQFT in short) ; see
\cite{Tu}.  Here $V_n$ is the $n$-dimensional irreducible
representation of $\mathfrak{sl}_2$ and $\xi$ is a complex root of
unity of order divisible by $4$. When $M=S^3$, and with the proper
normalization, we may identify the WRT invariant with the colored
Jones polynomial as follows:
$$
Z_{(S^3,K),\xi,V_n}= \ev_{\xi}J_K(n),
$$
where
$$
\ev_{\xi}: \BZ[q^{\pm 1/4}] \longto \BZ[\xi]
$$
is the {\em evaluation} $q^{1/4}=\xi$. From the point of view of
TQFT, it is not clear why the Jones polynomial is even a polynomial.

Perturbative quantum field theory also constructs an invariant
$\Zp_{(N,K),V_n}$ which takes values in the power series ring
$\BQ[[h]]$. When $M=S^3$, $\Zp_{(S^3,K),V_n}$ is the composition of
the {\em Kontsevich integral} of a knot with the $\mathfrak{sl}_2$
weight system, and when $N$ is arbitrary, $\Zp_{(N,K),V_n}$ is the
composition of the LMO invariant with the $\mathfrak{sl}_2$ weight
system, see \cite{LMO} and \cite{A}. When $M=S^3$, the main identity
is:
$$
\ev_h Z_K(n)=
\Zp_{(N,K),V_n} \in \BQ[[h]]
$$
where
$$
\ev_h : \BZ[q^{\pm 1/4}] \longto \BQ[[h]]
$$
is the evaluation $q=e^h$; see for example \cite{O}.
Thus, from the point of view of
perturbative TQFT the Jones polynomial is a formal power series in $h$.

\subsection{The Habiro ring}
\lbl{sub.habiro}

Despite the apparent failure of TQFT to explain the polynomial
aspect of the Jones polynomial, there is one gain. Namely let us fix
$(N,K)$ a natural number $n$ and let $\Om_4$ denote the set of
complex roots of unity whose order is divisible by $4$. Then, the
WRT invariant gives a function:
\begin{equation}
\lbl{eq.tau}
Z_{(N,K),V_n} : \Om_4 \longto \BC.
\end{equation}
This function is not continuous, and does not have nice analytic
properties.  Habiro introduced an alternative notion of  {\em
analytic functions}. The latter are by definition elements of the
{\em Habiro ring}, defined by:
\begin{equation}
\lbl{eq.BZhat}
\Ha := \lim_{\leftarrow n}\Z[q]/((1-q)(1-q^2)\dots(1-q^n)).
\end{equation}
The ring $\Ha$ can be considered as the set of all series of the
form
\begin{equation}
\lbl{eq.fq}
f(q) = \sum_{n=0}^\infty f_n(q) \, (1-q)(1-q^2)\dots(1-q^n),
\qquad \text{where } \quad f_n(q) \in \Z[q],
\end{equation}
with the warning that $f(q)$ does {\em not} uniquely determine $(f_n(q))$.

It turns out that the Habiro ring $\Ha$ shares many properties
with the ring of germs of complex analytic functions, (see \cite{Ha1})
and plays an important role in Quantum Topology.
In particular, elements of the Habiro ring
\begin{itemize}
\item[(a)]
can be differentiated with respect to $q$,
\item[(b)]
can be evaluated at the set $\Om$ of complex roots of unity,
\item[(c)]
have Taylor series expansions that uniquely determine them,
\item[(d)]
form an integral domain.
\end{itemize}
These properties suggest that we consider $\Ha$ as a class of ``analytic
functions" with domain $\Omega$. For proofs of these properties,
we refer the reader to \cite{Ha2}.

Let us comment a bit further on these properties. (a) is obvious from
Equation \eqref{eq.fq}. (b) also follows from \eqref{eq.fq} because
when $q$ is a root of unity,
only a finite number of terms in the right hand side of Equation
\eqref{eq.fq} are not $0$, hence $f(q)$ is defined as a complex number.
Thus one can consider every $f\in
\Ha$ as a function with domain $\Omega$ the set of roots of unity.
(a) and (b) imply that elements of the Habiro ring have Taylor series
expansions at every complex root of unity.
In particular, every $f\in \Ha$
has a Taylor expansion
$$
T_1(f) \in \BZ[[q-1]] \subset \BQ[[h]]
$$
(where $q=e^h$). What is nontrivial is the fact that $T_1(f)$ uniquely
uniquely determines $f$. (d) follows immediately from (c).
Another nontrivial fact is that
if $f(\xi)=g(\xi)$ at infinitely many roots
$\xi$ of prime power orders, then $f=g$ in $\Ha$. Therefore, one has the
following corollary of Habiro's results:

\begin{proposition}
\lbl{101}
The map $\Ha \to \prod_{n=1}^\infty \BZ[\xi_n]$, $f\to (f(\xi_1),
f(\xi_2),\dots)$ is injective. Here $\xi_n=\exp (2\pi i/n)$.
\end{proposition}

An example of a non-trivial element of the Habiro ring (of interest
to quantum topology) is the following element:
$$
f(q)=\sum_{n=0}^\infty  (1-q)(1-q^2)\dots(1-q^n)
$$
studied by Kontsevich (unpublished), Zagier and the second author;
see \cite{Za} and \cite{Le3}.

\subsection{Statement of the results}
\lbl{sub.results}

Ending our discussion on the Habiro ring, let us get back to
Equation \eqref{eq.tau}, and let us assume for simplicity that $K$
is a $0$-framed knot in an \ihs\ $N$. It turns out that there is a
function $f\in q^{n/2}\Ha$ such that for every $\xi\in \Omega_4$,

\begin{equation}
Z_{(N,K),V_n,\xi} = \ev_{\xi} f. 
\end{equation}

Proposition \ref{101} shows that this $f$ is an invariant of the
triple $N,K,n$, and we denote it by $J_{(N,K)}(n)$. Note that if $n$
is odd, then $J_{(N,K)}(n) \in \Ha$, otherwise $J_{(N,K)}(n) \in
q^{1/2}\Ha$. Our next definition answers  Question \ref{que.1} and
also explains the title of the paper:

\begin{definition}
\lbl{def.cjones}
The {\em colored Jones function} $J_{(N,K)}$ of a 0-framed knot $K$ in an
\ihs\ $N$ is defined by:
\begin{equation}
J_{(N,K)}: \BN \longto \Lhat,
\end{equation}
where
\begin{eqnarray}
\lbl{eq.Lhat}
\Lhat &:=& \Ha + q^{1/2}\,\Ha \\ \notag
& \cong & \Ha[x]/(x^2-q).
\end{eqnarray}
\end{definition}

\begin{remark}
\lbl{rem.framing}
If $K$ is an $0$-framed knot in an \ihs\ $N$
and $K_f$ denotes the knot $K$ with framing $f$ (where $f \in \BZ$),
then
$$
J_{(N,K_f)}=q^{f \frac{n^2-1}{4}} J_{(N,K)}.
$$
Therefore, the colored Jones function of a framed knot in an \ihs\ takes value
in the ring
$$
\Ha[x]/(x^4-q).
$$
\end{remark}

\subsection{The colored Jones function is $q$-holonomic}
\lbl{sub.qholo}

Now that we know what the Jones polynomial of a knot in an \ihs\ really
is, we may extend known results and conjectures
of the Jones polynomial of a knot in
$S^3$ to knots in \ihs s. One of these results, due to the authors,
is the fact that the colored Jones function of a knot in $S^3$ is
$q$-holonomic, see \cite{GL}. Let us recall this notion in our setting.

Define the linear operators $L$ and $M$ acting on maps $f$ from
$\BN$ to an $\BZ[q^{1/2}]$-module by:
$$
(M f)(n)= q^{n/2}f(n), \hspace{1cm} (L f)(n) = f(n+1).
$$
It is easy to see that $LM= q^{1/2} ML$, and that $L,M$ generate the
{\em quantum plane} $\cA$, the non-commutative ring with
presentation
$$
\cA=\BZ[q^{\pm 1/2}]\langle M,L \ra/(LM= q^{1/2} ML).
$$

The {\em recurrence ideal} of the discrete function $f$ is the left
ideal $I_f$ in $\cA$ that annihilates $f$:
$$
I_f=\{P \in \cA \, | \quad P\,f=0\}.
$$

We say that $f$ is $q$-{\em holonomic}, or $f$ satisfies a linear
recurrence relation, if $I_f\neq 0$. In \cite{GL} we proved that for
every knot $K$ in $S^3$, the function $J_K$ is $q$-holonomic.
In other words, $J_K$ satisfies a linear recursion relation with
coefficients Laurent polynomials  in $q^{1/2}$ and $q^{n/2}$.

\begin{theorem}
\lbl{thm.1} For every 0-framed knot $K$ in an \ihs\ $N$, the colored
Jones function $J_{(N,K)}$ is $q$-holonomic.
\end{theorem}

The main ideas behind the proof of the above theorem is that
\begin{itemize}
\item
every pair $(N,K)$ as above is obtained from unit-framed surgery
from an algebraically split link $K \cup L$ in $S^3$.

\item
The function $J_{(N,K)}$ is obtained from the colored Jones function
$J_{K \cup L}$ by elimination of the variables corresponding to $L$.
\item
The function $J_{K \cup L}$ is $q$-holonomic in all its variables
(see section \ref{holonomic}  below), by \cite{GL}.
\item
Elimination preserves $q$-holonomicity.
\end{itemize}

\subsection{The recurrence polynomial}
\lbl{sub.recurrence}

Let $I_{(N,K)}$ denote the recurrence ideal of $J_{(N,K)}$. Then
$I_{(N,K)}$ is a left ideal of  $\cA$, which is not a principal
ideal domain. Hence $I_{(N,K)}$ might not be generated by a single
element. The first author \cite{Ga} noticed that by adding to $\cA$
all the inverses of polynomials in $M$ one gets a principal ideal
domain $\tilde\cA$, and hence from the ideal $I_{(N,K)}$ one can
define a polynomial invariant. Formally, let $\BQ(q^{1/2},M)$ be the
fractional field of the polynomial ring $\cR[M]$. Let $\tilde \cA$
be the set of all polynomials in the variable $L$ with coefficients
in $\BQ(q^{1/2},M)$:
$$
\tilde\cA =\{\sum_{k=0}^\infty a_k(M) L^k \,\, | \quad a_k(M)\in
\BQ(q^{1/2},M),
\,\,\, a_k=0  \quad k \gg 0 \},
$$
with multiplication given by 
$$
a(M) L^k \cdot b(M)L^l=a(M)\, b(q^{k/2}M) L^{k+l}.
$$
It is known that $\tilde\cA$ is a twisted polynomial ring (an
{\em Ore extension} of
$\BQ(q^{1/2},M)$),
and consequently $\tilde\cA$ is a principal left-ideal domain,
and $\cA$ embeds as a subring of $\tilde \cA$. The ideal extension
$\tilde I_{(N,K)}  := \tilde \cA \, I_{(N,K)}$  is then generated by
a single polynomial

$$\alpha_{(N,K)}(L,M,q) = \sum_{i=0}^n \alpha_{(N,K),i}(M,q) L^i$$
where the degree in $L$ is assumed to be minimal and all the
coefficients $\alpha_{(N,K),i}(M,q)\in \BZ[q^{\pm1/2},M]$ are
assumed to be co-prime. That $\alpha_{(N,K)}$ can be chosen to have
integer coefficients follows from the fact that $J_{(N,K)}(n) \in
\Lhat$. It is clear that $\alpha_{(N,K)}(L,M,q)$ annihilates
$J_{(N,K)}$, and hence it is in the recurrence ideal $I_{(N,K)}$.
Note that $\alpha_{(N,K)}(M,L;q)$ is defined up to a factor $\pm
q^{a/2} M^{b}, a,b\in \BZ$. We will call $\alpha_{(N,K)}$ the {\em
recurrence polynomial} of $(N,K)$.

Let us say that two non-zero polynomials $p_1$ and $p_2$ in
variables $L$ and $M$ are {\em $M$-essentially equal} (and write
$p_1 \,\overset{M}{=}\, p_2$) if their ratio is a function of $M$
alone. The next conjecture generalizes the AJ Conjecture of the
first author (see \cite{Ga}):

\begin{conjecture}
\lbl{conj.AJ}
For every knot $K$ in an \ihs \ $N$, we have:
$$
\alpha_{(N,K)}(L,M,1) \,\overset{M}{=}\ A_{(N,K)}(L,M),
$$
where $A_{(N,K)}(L,M)$ is the $A$-polynomial of $(N,K)$ defined
by \cite{CCGLS}.
\end{conjecture}

For some partial results confirming the conjecture, see
\cite{Hi,Le2,Ga}. In particular, in \cite{Le2} the second author  used
Kauffman bracket modules to established the conjecture for a large
class of 2-bridge knots in $S^3$. The above conjecture compares a
quantum invariant (the recurrence ideal of the colored Jones
function) with a classical one (namely the $A$-polynomial). One
motivation of the conjecture is the dream of {\em quantization} and
{\em semiclassical analysis}, in the context of 3-manifolds with
torus boundary. Another motivation is the fact that the Kauffman
bracket skein module is in a sense a quantization of the coordinate
ring of the $\SL_2(\BC)$ character variety of a 3-manifold; see
\cite{PS}. The $A$-polynomial is the coordinate ring of the
$\SL_2(\BC)$ character variety of the knot complement, viewed from
the boundary torus. On the other hand, the recurrence polynomial is
in a sense a quantization of the classical coordinate ring. Thus, we
are back to the Kauffman bracket skein module. And with this happy
thought in mind, we end this section.

\section{A review of Habiro's work}
\lbl{sec.habiro}

In this section we review Habiro's work on the integrality of the
colored Jones polynomial and the WRT invariants of links in \ihs s.
For a detailed discussion, we refer the reader to \cite{Ha1} and
\cite{Ha3}.

\subsection{The quantum group $U_q(\fsl_2)$}
\lbl{sub.qgroups}

We begin by reviewing some necessary representation theory of
quantum groups.

We will use the following analogs of {\em quantum integers} and {\em quantum
factorials}:
$$
\{n\}=q^{n/2}-q^{-n/2}, \qquad \{n\}!=\prod_{i=1}^n \{i\}! .
$$
Let
\begin{eqnarray*}
C(n,k) &:=& \prod_{j=n-k}^{n+k} \{j\} = \frac{\{n+k\}!}{\{n-k-1\}!}
\end{eqnarray*}

Notice that $C(n,k)=0$ for $k>n$. Also $C(n,k) \in q^{n/2}\BZ[q^{\pm
1}]$ and hence $\{n\}C(n,k) \in \BZ[q^{\pm 1}]$.

Consider the quantized enveloping
algebra $U_q(\fsl_2)$ be
defined as Jantzen's book \cite{Jantzen}, except that our $q$ is
equal to $q^2$ of \cite{Jantzen}, and our ground field is extended
to $\BC(q^{1/4})$, instead of $\BC(q^{1/2})$ as in \cite{Jantzen}.

The theory of type 1 $U_q(\fsl_2)$-modules is totally similar to that
of $\fsl_2$. For each positive integer there is a unique irreducible
type 1 $U_q(\fsl_2)$-module $V_n$, , see eg \cite{Jantzen} (where
$V_n$ is denoted by $L(n,+)$). Let $\RR$ be the Grothendieck ring of
finite-dimensional type 1 $U_q(\fsl_2)$-modules, tensored by $
\BC(q^{1/4})$. Then $\RR$ is freely spanned by $V_1,V_2,\dots$, and
as an algebra is isomorphic to the polynomial algebra
$\BC(q^{1/4})[V_2]$. Let $\RR^{e}$ denote the subspace of $\RR$
spanned by even powers of $V_2$. For every framed, $r$-component,
oriented link $\cL$ whose components are ordered, the quantum
invariant $J_\LL$, extended linearly, can be considered as a
function from $\RR^r$ to $\BC(q^{1/4})$.

There is a symmetric bilinear form on $\RR$:
\begin{equation}
\lbl{eq.pairing} \RR \times \RR \longto \BC(q^{1/4}), \qquad \langle
V, U \ra := \psdraw{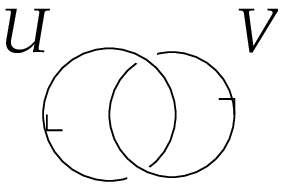}{0.40in} .
\end{equation}
Here the right hand side is the colored Jones polynomial of the Hopf
link, colored by $V$ and $U$. For example, one has $\la V_n, V_m\ra
= [nm]$.

Habiro defines the following elements in $\RR$:
$$
S_k := \prod_{i=1}^k(V_2^2- \{i\}^2) , \quad \text{and} $$
$$ P_k :=
\prod_{i=1}^{k}(V_2-q^{(2i-1)/2}-q^{-(2i-1)/2}), \quad P_k':=
\frac{P_k}{\{k\}!}, \quad P_k'':= \frac{\{1\}}{\{2k+1\}!} P_k.
$$
It is clear that $\{ S_k \}_{k \in \BN}$ and $\{ P_k \}_{k \in \BN}$
are bases of $\RR^e$ and $\RR$ respectively. Moreover, the two bases
are dual under the pairing \eqref{eq.pairing}. That is, we have:
\begin{eqnarray}
\lbl{eq.dualb}
\langle P''_k, S_n \ra &=& \d_{n,k} .
\end{eqnarray}

It's easy to show  that
\begin{equation}
\lbl{106}
\langle V_n, S_k \ra = \frac{C(n,k)}{\{1\}} .
\end{equation}

\subsection{Universal quantum invariant of  knots}
\lbl{sub.universalqknots}

Suppose $K$ is a 0-framed oriented knot in $S^3$. Then $J_K$ can be
considered as a $\BC(q^{1/4})$-linear map from $\RR$ to the ground
field $\BC(q^{1/4})$. Using the bilinear form \eqref{eq.pairing} one
can also consider each $S_k$ as a linear form on $\RR$. Habiro
showed that $J_K$ is always an infinite linear combination of the
$S_k$, i.e. for every $V\in \RR$,

\begin{equation}
\lbl{eq.basica} J_K(V)=\sum_{k=0}^\infty C_K(k) \langle V, S_k \ra.
\end{equation}

It is easy to show that for every $V\in \RR$, if $k$ is big enough
then $\la V,S_k\ra =0$, hence the right hand side has a meaning for
every $V\in \RR$. Using the orthogonality \eqref{eq.dualb} one see
that

$$
 C_K(k) = J_K( P''_k)
$$

A priori $C(k)$ belongs to the ground field $\BC(q^{1/4})$. A
difficult result of Habiro is that $C(k)$ is always a Laurent
polynomial in $q$ with integer coefficients, $C(k) \in \BZ[q^{\pm
1}]$.

Using Equation \eqref{eq.basica}  for $V=V_n$ and \eqref{106}, one
has
\begin{equation}
\lbl{eq.JK} J_K(n)= \sum_{k=0}^\infty C_K(k) C(n,k)/\{1\}.
\end{equation}

\subsection{The case of links}
\lbl{sub.links}

Suppose $\LL$ is an algebraically
split $r$-component link in $S^3$ with 0 framing on each component.
Then it is known that $J_{\LL}(n_1,\dots,n_r)$ is in
$\BZ[q^{\pm1}]$ or in $q^{1/2}\BZ[q^{\pm1}]$, according as
$n_1+\dots+n_r-r$ is even or odd.

Equation (\ref{eq.basica}) can be generalized to the link case,
and we have the fundamental equation:
\begin{equation}
\lbl{eq.basiclink} J_\LL(W_1,\dots,W_r)
 = \sum_{k_i=0}^\infty C_\LL(k_1,\dots,k_r)
\langle W_1, S_{k_1}\ra \dots \langle W_r, S_{k_r}\ra  ,
\end{equation}
for every $W_1,\dots,\W_r\in \RR$, and
$$
C_\LL(k_1,\dots,k_r)=J_\LL(P''_{k_1}, \dots , P''_{k_r}).
$$
The analogue of Equation \eqref{eq.JK} is:
\begin{eqnarray}
\lbl{eq.JL} J_{\LL}(n_1,\dots,n_r) &=& \sum_{k_1,\dots,k_r=0}^\infty
C_{\LL}(k_1,\dots,k_r) \prod_{i=1}^r \frac{C(n_i,k_i)}{\{1\}}.
\end{eqnarray}

Unlike the knot case, the coefficient $C_{\LL}(k_1,\dots,k_r)$ of
the right hand side is not always a Laurent polynomial, but rather a
rational function in $q$. One of the main results of Habiro is the
integrality of a minor modification $\ti C_L$ of the cyclotomic
function $C_L$. Let us define

\begin{eqnarray*}
\ti C_\LL(k_1,\dots,k_r) &:=& J_\LL(P'_{k_1}, \dots , P'_{k_r}) \\
 &=&
C_\LL(k_1,\dots,k_r) \frac{\{2k_1+1\}!}{\{k_1\}!\{1\}} \dots
\frac{\{2k_r+1\}!}{\{k_r\}!\{1\}},
\end{eqnarray*}
 Habiro proves that

\begin{theorem}\cite[Thm.3.3]{Ha1}
\lbl{thm.habiro1}
If $L$ is algebraically split
and zero framed link in $S^3$, then
$$
\ti C_L (k_1,\dots,k_r) \in
\frac{\{2m+1\}!}{\{m\}!\{1\}} \,\,\BZ[q^{\pm1/2}],
$$
where $m=\max\{k_1,\dots, k_r\}$.
\end{theorem}

The important thing is that for every $n$,
$\frac{\{2n+1\}!}{\{n\}!\{1\}}$ is divisible by $(1-q) \dots
(1-q^n)$. This guarantees that for any sequence $f(k) \in
\BZ[q^{\pm1}]$ (for $k \in \BN$), the series
$$
\sum_{k=0}^\infty f(k) \frac{\{2k+1\}!}{\{k\}!\{1\}}
$$
converges in $\Lhat$.

\subsection{Invariants of \ihs s}
\lbl{sub.ihs}

Suppose $N$ is an \ihs, which   is obtained by surgery on $S^3$
along an algebraically split $r$-component link $\LL$, with framing
$f_1=\pm1,\dots,f_r=\pm1$. Let $\LL^{(0)}$ be the link $\LL$ with
all framing switched to 0.

 Habiro introduced
the following elements $\w^{\pm 1}$ in some completion of $\RR$:

\begin{eqnarray*}
\w &=& \sum_{j=0}^\infty q^{j(j+3)/4} P_j' \\
\w^{-1} &=& \sum_{j=0}^\infty (-1)^j q^{-j(j+3)/4} P_j'
\end{eqnarray*}

Although each of $w^{\pm 1}$ is an infinite sum of elements in
$\RR$, Theorem \ref{thm.habiro1} (with the remark after the theorem)
shows that $J_{\LL^{(0)}}(w^{\pm 1}, \dots, w^{\pm 1})$ always
belongs to $\Ha$.

Using a version of Kirby's calculus for algebraically split links,
Habiro proved that
$$
J_N:=J_{\LL^{(0)}}(w^{-f_1}, \dots, w^{-f_r})\in \Ha
$$
is an invariant of the \ihs\ $N$, i.e. does not depend on the
choice of the surgery link $\LL$. Moreover, one has the following

\begin{proposition} { \rm (Habiro)}
 The evaluation of $J_N$ (as an element of $\Ha$) at a root of unity
coincides with the WRT invariant $Z_N$ at that same root of unity.
\label{103}
\end{proposition}
\begin{proof}
We give a proof for this fact here, since we will adapt the proof
for the relative case later.

Orthogonality (that is Equations \eqref{eq.dualb} and \eqref{106})
implies that
$$
\langle P'_k, S_n \ra = \d_{n,k} \frac{\{2k +1\}!}{\{k \}!\{1\}}.
$$
Thus, for $f=\pm 1$ we have:
\begin{eqnarray*}
\langle \w^{-f}, S_k \ra = \langle (-f q)^{-f k(k+3)/4} P'_k, S_k
\ra  = (-f q)^{-f k(k+3)/4} \frac{\{2k +1\}!}{\{k \}!\{1\}}
\end{eqnarray*}
This, together with Equation \eqref{eq.basiclink} implies that
\begin{eqnarray*}
J_N & = & \sum_{k_1,\dots, k_r=0}^\infty  C(k_1,\dots, k_r) \langle
\w^{-f_1}, S_{k_1}
\ra \dots \langle \w^{-f_r}, S_{k_r} \ra \\
&=& \sum_{k_1,\dots, k_r=0}^\infty  C(k_1,\dots, k_r)
\prod_{i=1}^r\frac{\{2k_i +1\}!}{\{k_i \}!\{1\}}  (-f_i q)^{-f_i
k_i(k_3+3)/4} \end{eqnarray*}

Hence
\begin{equation}
 J_N  = \sum_{k_1,\dots, k_r=0}^\infty \ti
C(k_1,\dots, k_r) \prod_{i=1}^r (-f_i q)^{-f_i k_i(k_3+3)/4}.
\label{108}
\end{equation}

On the other hand, if $q^{1/4}$ is a root of unity of order $4d$,
then the quantum invariant $Z_{N}(q)$ is given by (see eg \cite{KM})

$$ Z_{N}(q) = \frac{\sum_{n_i=1}^{4d} J_{\cL^{(0)}}(n_1,\dots,n_r) \,
\prod_{i=1}^r [n_i] q^{\frac{f_i}{4}(n_i^2-1)}} {\prod_{i=1}^{r}
\left (  \sum_{j=1}^{4d} [n_i]^2 q^{\frac{f_i}{4}(n_i^2-1)}
\right)},$$ which, by \eqref{eq.JL} and Corollary \ref{105} below,
is equal to the right hand side of \eqref{108}.
\end{proof}

\subsection{Laplace transform} Suppose $q^{1/4}$ is a primitive
root of unity of order $4d$. For $f=\pm 1$, Let $\gamma_f$ be the
Gauss sum
$$\gamma_f:= \sum_{k=1}^{4d} q^{\frac{f}{4}(k^2-1)}.$$
 The exact value of $\gamma_f$ can be calculated, see eg. \cite{KM},
but for us
 it is only important that $\gamma_f \neq 0$.

We define the Laplace transform $\Lp_{f;n}$ as follows

$$ \Lp_{f;n}(q^{na +b}) := q^{b-fa^2}.$$

\begin{lemma}
For any integers $a,b$ we have
$$
\sum_{n=1}^{4d} q^{\frac{f}{4}(n^2-1)}\, q^{b+na} =
\gamma_f\, \Lp_{f;n}(q^{b+na}).
$$
\end{lemma}

\begin{proof}
We have
$$
q^{\frac{f}{4}(n^2-1)} \, q^{b+na} = q^{\frac{f}{4}((n
+2a)^2-1)} \, q^{b-fa^2}.
$$
Summing up $n$ from 1 to $4d$ we get the result.
\end{proof}

The following lemma, in another notation, was Lemma 2.2 of \cite{BBL}:

\begin{lemma}
One has
$$
\Lp_{f,n}([n] C(n,k)/\{1\}) =
2(q^{-f}-1) (-fq)^{-fk(k+3)/4}\,\frac{\{2k+1\}!}{\{k\}! \{1\}}.
$$
\end{lemma}

\begin{corollary}
For $q^{1/4}$ a root of unity of order $4d$, one has
$$
\lbl{105}
\frac{\sum_{n=1}^{4d} q^{\frac{f}{4}(n^2-1)}\, [n]
C(n,k)/\{1\}}{\sum_{n=1}^{4d} q^{\frac{f}{4}(n^2-1)}\, [n]^2}
= (-fq)^{-fk(k+3)/4}\,\frac{\{2k+1\}!}{\{k\}!
\{1\}}.
$$
\end{corollary}

\begin{remark}
\lbl{rem.altformula}
Despite the appearance of powers of $q^{1/4}$ in $\w^{\pm}$, $J_N$ contains
integral powers of $q$. This follows from
\begin{eqnarray}
\lbl{eq.alt1}
q^{n(n+3)/4} \frac{\{2n+1\}!}{\{n\}!\{1\}} &=& (-1)^n q^{-n(n+1)/2}
\frac{\{2n+1\}_{-} !}{\{n\}_{-} !\{1\}_{-}} \\
\lbl{eq.alt2}
(-1)^n q^{-n(n+3)/4} \frac{\{2n+1\}!}{\{n\}!\{1\}} &=& q^{-n(n+2)}
\frac{\{2n+1\}_{-} !}{\{n\}_{-} !\{1\}_{-}} ,
\end{eqnarray}
where the {\em unbalanced} quantum integers and factorials are defined by:
$$
\{n \}_{-}=1-q^n, \qquad \{n \}_{-} ! = \{1 \}_{-} \dots \{n \}_{-} .
$$
\end{remark}

The next corollary gives a formula for $J_N$, where $N$ is obtained by
$\pm 1$ surgery on a knot $K$ in $S^3$:

\begin{corollary}
\lbl{cor.pmsurgery}
If $S^3_{K,\pm 1}$ denotes the result of $\pm 1$ surgery on a knot $K$
in $S^3$, then
\begin{eqnarray}
\lbl{eq.altm}
J_{S^3_{K,-1}} &=& \sum_{k=0}^\infty J_K(P''_k)
q^{k(k+3)/4} \frac{\{2k+1\}!}{\{k\}!\{1\}} \\
\notag
&=& \sum_{k=0}^\infty J_K(P''_k)
(-1)^k q^{-k(k+1)/2}
\frac{\{2k+1\}_{-} !}{\{k\}_{-} !\{1\}_{-}} \\
\lbl{eq.altp}
J_{S^3_{K,+1}} &=& \sum_{k=0}^\infty J_K(P''_k)
(-1)^k q^{-k(k+3)/4} \frac{\{2k+1\}!}{\{k\}!\{1\}} \\
\notag
&=& \sum_{k=0}^\infty J_K(P''_k)
q^{-k(k+2)}
\frac{\{2k+1\}_{-} !}{\{k\}_{-} !\{1\}_{-}}.
\end{eqnarray}
\end{corollary}

\section{A relative version}
\lbl{sec.relative}

Suppose  $\LL$ is an algebraically split $r$-component link in an
\ihs\ $N$, each component of which  has framing 0. It is known that
there is an algebraically split framed link in $S^3$, which is the
disjoint union of 2 sublinks $ \cL_1$ and $\cL_2$, such that the
framing of each component of $\cL_2$ is $\pm 1$, and that the
surgery along $\cL_2$ transforms $(S^3,\cL_1)$ to $(N,\cL)$. Then
$\cL_1$ has $r$ components, each with framing 0. Assume that $\cL_2$
has $s$ components whose framings are $f_1, f_2,\dots, f_s$. Let
$\cL^{(0)}$ be the link $\cL_1\cup \cL_2$ with all the framings
switched to 0.

Suppose $n_1,\dots,n_r$ are positive integer. Let $Z_{N,
\cL}(n_1,\dots,n_r;\xi)$ be the WRT invariant of the pair $(N,\LL)$,
when $\LL$ is colored by the $U_q(\fsl_2)$-module
$V_{n_1},\dots,V_{n_r}$, at the root $\xi$ of unity of order
divisible by 4. Our $\xi$ here (which is $q^{1/4}$)  is equal to $t$
in \cite{KM}.

\subsection{WRT invariant $Z_{(N,\cL)}$ as element of $\Lhat$}

\begin{lemma} For  positive integer $n_1,\dots,n_r$, one has
$J_{\cL^{(0)}}(n_1,\dots,n_r,\w^{-f_1}, \dots,
\w^{-f_s})$ belongs to $\Lhat$.
\end{lemma}

\begin{proof}
Using (\ref{eq.basiclink}) and a  computation analogous to the one
in the proof of Proposition \ref{103}, one has
\begin{equation}
J_{\cL^{(0)}}(n_1,\dots,n_r,\w^{-f_1}, \dots, \w^{-f_s})=
 \sum_{0\le k_i < n_i} \left( \prod_{i=1}^r \frac{
\,C(n_i,k_i)}{\{1\}}\right) \, D(k_1,\dots,k_i) \label{100}
\end{equation}
where
\begin{eqnarray}
\lbl{eq.compute12} D(k_1,\dots,k_r) &=&
\sum_{l_1,\dots,l_s=0}^\infty \left( \prod_{i=1}^s (-f_i q)^{-f_i
l_i(l_i+3)/4}\right) \, J_{\cL^{(0)}}(P''_{k_1},\dots, P''_{k_r},
P'_{l_1}, \dots, P'_{l_s}).
\end{eqnarray}

Theorem \ref{thm.habiro1} implies that
\begin{equation}
\lbl{eq.integrality} J_{\cL^{(0)}}(P''_{k_1},\dots, P''_{k_r},
P'_{l_1}, \dots, P'_{l_s}) \in \left( \prod_{i=1}^r
\frac{\{k_i\}!\{1\}}{\{2k_i+1\}!}\right)\,
\frac{\{2m+1\}!}{\{m\}!\{1\}} \,\,\cR
\end{equation}
where $m=\max\{k_1,\dots, k_r, l_1,\dots,l_s\}$. Using the easily
verified identity

$$  \frac{C(n,k)}{\{1\}} \frac{\{k\}!\{1\}}{\{2k+1\}!}
= \qbinom{n+k}{2k+1}\, \{k\}!  $$

we see that

$$ \prod_{i=1}^r \frac{
\,C(n_i,k_i)} {\{1\}} \, J_{\cL^{(0)}}(P''_{k_1},\dots, P''_{k_r},
P'_{l_1}, \dots, P'_{l_s}) \in \left( \prod_{i=1}^r\{k_i\}!\right)
\frac{\{2m+1\}!}{\{m\}!\{1\}} \,\,\cR.$$

Hence for fixed $k_1,\dots,k_r$, the term in the sum of the right
hand side of (\ref{100}) is in $\Lhat$.
\end{proof}

Recall that if $q^{1/4}$ is a root of unity of order $4d$, then the
WRT
 invariant is defined by

$$ Z_{(N,\cL)}(n_1,\dots,n_r; q) = \frac{\sum_{n_i=1}^{4d}
J_{\cL^{(0)}}(n_1,\dots,n_r) \,
\prod_{i=1}^r [n_i] q^{\frac{f_i}{4}(n_i^2-1)}} {\prod_{i=1}^{r}
\left (  \sum_{j=1}^{4d} [n_i]^2 q^{\frac{f_i}{4}(n_i^2-1)}
\right)}.$$

Hence the proof of Proposition \ref{103} can be easily generalized
to the relative case, and one gets

\begin{proposition}
\lbl{prop.eval}
The evaluation of $J_{\cL^{(0)}}(n_1,\dots,n_r,\w^{-f_1}, \dots,
\w^{-f_s})$ at a root $q^{1/4}=\xi$ of unity coincides with the
quantum invariant $Z_{(N,\cL)}(n_1,\dots,n_r)$ at that same root of
unity. \newline
Hence $J_{\cL^{(0)}}(n_1,\dots,n_r,\w^{-f_1}, \dots,
\w^{-f_s})$ is an invariant of the link $\cL$ in $N$, colored by
$n_1,\dots,n_r$.
\end{proposition}

\begin{remark}
Habiro's argument can directly show that
$J_{\cL^{(0)}}(n_1,\dots,n_r,\w^{-f_1}, \dots, \w^{-f_s})$ is an
invariant of the link $\cL$ in $N$, colored by $n_1,\dots,n_r$.
\end{remark}

Let us denote $J_{\cL^{(0)}}(n_1,\dots,n_r,\w^{-f_1}, \dots,
\w^{-f_s})$ by $J_{(N,\cL)}(n_1,\dots,n_r)$, which is an element of
$\Lhat$. We can consider $J_{(N,\cL)}$ as a function from $\BN^r$ to
$\Lhat$.

\section{$q$-holonomicity in many variables}
\lbl{sec.qholomany}

Theorem \ref{thm.1} is a special case of Theorem \ref{109} below,
which follows from the fact that the quantum invariants can be built
from elementary blocks that are $q$-holonomic, and the operations
that patch the blocks together to give the colored Jones function
preserve $q$-holonomicity.  First we need the notion of
$q$-holonomicity in many variables, introduced by Sabbah \cite{C},
generalizing Bernstein's notion of (usual) holonomicity
\cite{B1,B2}.

\subsection{$q$-holonomicity in many variables}  Consider the
operators $L_i$ and $M_j$ for $1 \leq i,j \leq r$ which act on
functions $f$ from $\BN^r$ to a $\BZ[q^{\pm 1/2}]$-module by
\begin{eqnarray*}
(M_i f)(n_1,\dots,n_r)&=&q^{n_i/2}f(n_1,\dots, n_r) \\
(L_i f)(n_1,\dots,n_r) &=& f(n_1,\dots, n_{i-1}, n_i+1, n_{i+1},
\dots, n_r).
\end{eqnarray*}
It is easy to see that the following relations hold:
\begin{equation*}
\tag{$\text{Rel}_{q}$}
\begin{aligned}
M_i M_j&= M_j M_i  & L_i L_j &= L_j L_i \\
M_i L_j&= L_j M_i \, \text{for} \, i \neq j & L_i M_i &= q^{1/2} M_i
L_i
\end{aligned}
\end{equation*}
We define the $r$-dimensional quantum space $\A_r$ to be a
noncommutative algebra with presentation
$$
\A_r= \frac{\BZ[q^{\pm1/2}]\la M_1, \dots, M_r , L_1, \dots, L_r
\ra}{ (\text{Rel}_q)}.
$$

For a function $f$ as above one can define the left ideal $\calI_f$
in $\A_r$ by
$$
\calI_f := \{P \in \A_r | Pf=0 \}    .
$$
If we want to determine
a function $f$ by a finite list of
initial conditions, it does not suffice to ensure that $f$ satisfies
one nontrivial recursion relation if $r\ge 2$. The key notion that
we need instead is $q$-holonomicity. Intuitively, $f$ is
$q$-holonomic if it satisfies a {\em maximally overdetermined
system} of linear difference equations with polynomial coefficients.
The exact definition of holonomicity is through homological
dimension, as follows.

Suppose  $I$ is a left $\A_r$-module. Let $F_m$ be the sub-space of
$\A_r$ spanned by polynomials in $M_i,L_i$ of total degree $\le m$.
Then the module $\A_r/I$ can be approximated by the sequence
$F_m/(F_m\cap I), m=1,2,...$. It turns out that, for $m
>>1$,  the dimension (over the fractional field $\BQ(q^{1/2})$) of
$F_m/(F_m\cap I)$ is a polynomial in $m$ whose degree is called the
{\em homological dimension} of $\A_r/I$ and is denoted by
$d(\A_r/I)$.

Bernstein's {\em famous inequality} (proved by Sabbah in the
$q$-case) states that the dimension of a non-0 module is $ \geq r$,
if the module has {\em no monomial torsions}, i.e., any non-trivial
element of the module  cannot be annihilated by a monomial in
$M_i,L_i$. Note that the left $\A_r$-module $ \A_r/\calI_f$ does not
have monomial torsion.

We say that a discrete function $f$ is
$q$-holonomic if $d(\A_r/\calI_f)\le r$. Note that if $f$ is
$q$-holonomic,  then by Bernstein's inequality, either
$\A_r/\calI_f=0$ or $d(\A_r/\calI_f)=r$. The former can happen only
if $f=0$.

\subsection{Assembling $q$-holonomic functions} \lbl{sub.assemble}

Here are some important operations that preserve $q$-holonomicity:

\begin{itemize}
\item
Sums and products of $q$-holonomic functions are $q$-holonomic.
\item
Specializations and extensions of $q$-holonomic functions are
$q$-holonomic. In other words, if $f(n_1,\dots,n_m)$ is
$q$-holonomic, the  so are the functions $g(n_2,\dots,n_m):=
f(a,n_2,\dots,n_m)$ and $h(n_1,\dots,n_m,n_{m+1}):=
f(n_1,\dots,n_m).$
\item
Diagonals of $q$-holonomic functions are $q$-holonomic. In other
words, if $f(n_1,\dots,n_m)$ is $q$-holonomic, then  so is the
function
$$g(n_2,\dots,n_m) := f(n_2,n_2,n_3,\dots,n_m).$$
\item
Linear substitution. If $f(n_1,\dots,n_m)$ is $q$-holonomic, then so
is the function, $g(n_1',\dots,n'_{m'})$, where each $n'_j$ is a
linear function of the $n_i$'s.
\item
Multisums of $q$-holonomic functions are $q$-holonomic. In other
words, if $f(n_1,\dots,n_m)$ is $q$-holonomic, the  so are the
functions $g$ and $h$, defined by
$$
g(a,b,n_2,\dots,n_m):= \sum_{n_1=a}^b f(n_1,n_2,\dots,n_m), \qquad
h(a,n_2,\dots,n_m) := \sum_{n_1=a}^\infty f(n_1,n_2,\dots,n_m)
$$
(assuming that the latter sum is finite for each $a$).
\end{itemize}

For a user-friendly explanation of these facts and for many
examples, see \cite{Ze} and \cite{PWZ}.

\subsection{Examples of $q$-holonomic functions}

The following functions are $q$-holonomic:
$$
n \to \{n\}, \qquad  n \to [n], \qquad   n\to [n]! :=\prod_{i=1}^n
[i], \qquad  n\to \{n\}!
$$
$$
(n,k) \to \{n\}_k := \begin{cases}
\prod_{i=1}^k \{n-i+1\},  & \text{if $k\ge 0$}\\
0 & \text{if $k<0$} \end{cases}
$$
$$
(n,k) \to \qbinom{n}{k} := \begin{cases}
\frac{\{n\}_k}{\{k\}_k}  & \text{if $k\ge 0$}\\
0 & \text{if $k<0$} \end{cases}.
$$

Also  $q$-holonomic is the delta function $\delta_{n,k}$. In fact,
we will encounter only sums, products, extensions, specializations,
diagonals, and multisums of the above functions.

\subsection{$q$-holonomicity of quantum invariants}
\lbl{holonomic}

\begin{theorem}
\lbl{109}
For a $0$-framed, algebraically split, $r$-component,
oriented link $\cL$ in an
\ihs\ $N$, the function $J_{(N,\cL)}: \BZ^r \to \Lhat$ is
$q$-holonomic.
\end{theorem}

\begin{proof}
From \cite{GL} we know that the function
$$
R(n,k) := (-1)^{n+1-k} \frac{\{1\} \{2k\}}{\{n+1-k\}!
\{n+1+k\}!}
$$
is $q$-holonomic, and  that
$$
P''_n = \sum_{k=1}^{n+1} R(n,k) V_k.
$$

By the result of \cite{GL} we know that
$J_{\cL^{(0)}}(k_1,\dots,k_r,l_1,\dots,l_s)$ is $q$-holonomic in all
variables,
hence $ J_{\cL^{(0)}}(P''_{k_1},\dots, P''_{k_r}, P'_{l_1}, \dots,
P'_{l_s})$ is $q$-holonomic in all variables
$k_1,\dots,k_r,l_1,\dots,l_s$. It follows from (\ref{eq.compute12})
that $D(k_1,\dots,k_r)$ is $q$-holonomic, since $(-f_1 q)^{-f_1
k_1(k_1+3)/4} \dots (-f_r q)^{-f_r k_r(k_r+3)/4}$ is $q$-holonomic
in all variables. Then equation (\ref{100}) shows that $J_{(N,\cL)}$
is $q$-holonomic.
\end{proof}

\ifx\undefined\bysame
    \newcommand{\bysame}{\leavevmode\hbox
to3em{\hrulefill}\,}
\fi

\end{document}